\documentclass[leqno, 12pt]{article}
\usepackage{amsmath,amsfonts,amsthm,amssymb,indentfirst}

\newtheorem{theorem}{Theorem}
\newtheorem{lemma}[theorem]{Lemma}
\newtheorem{definition}[theorem]{Definition}
\newtheorem{corollary}[theorem]{Corollary}
\newtheorem{proposition}[theorem]{Proposition}

\newtheorem{question}[theorem]{Question}

\theoremstyle{definition}

\newcommand{\Sym}{\mathrm{Sym}}
\newcommand{\Aut}{\mathrm{Aut}}
\newcommand{\Self}{\mathrm{Self}}

\newcommand{\Q}{\mathbb{Q}}
\newcommand{\Z}{\mathbb{Z}}

\begin{document}

\title{Groups where free subgroups are abundant}

\author{Zachary Mesyan\thanks{This work was done while the author was supported by a Postdoctoral Fellowship from the Center for Advanced Studies in Mathematics at Ben Gurion University, a Vatat Fellowship from the Israeli Council for Higher Education, and ISF grant 888/07.}}

\maketitle

\begin{abstract}
Given an infinite topological group $G$ and a cardinal $\kappa > 0$, we say that $G$ is \emph{almost $\kappa$-free} if the set of $\kappa$-tuples $(g_i)_{i \in \kappa} \in G^\kappa$ which freely generate free subgroups of $G$ is dense in $G^\kappa$. In this note we examine groups having this property and construct examples. For instance, we show that if $G$ is a non-discrete Hausdorff topological group that contains a dense free subgroup of rank $\kappa > 0$, then $G$ is almost $\kappa$-free. A consequence of this is that for any infinite set $\Omega$, the group of all permutations of $\Omega$ is almost $2^{|\Omega|}$-free. We also show that an infinite topological group is almost $\aleph_0$-free if and only if it is almost $n$-free for each positive integer $n.$ This generalizes the work of Dixon and Gartside-Knight.
\end{abstract}

\section{Introduction}

In 1990 Dixon~\cite{Dixon} showed that almost all finite sequences of permutations of a countably infinite set freely generate free subgroups. More precisely, letting $S = \Sym(\Z_+)$ denote the group of all permutations of the positive integers, Dixon showed that for each integer $n \geq 2$, the set $\{(g_1, \dots, g_n) \in S^n : \{g_1, \dots, g_n\} \ \text{freely generates a free subgroup of} \ S\}$ is comeagre in $S^n$, where $S$ is viewed as a topological group under the function topology. (The notion of ``comeagre" will be recalled in Section~\ref{generalities}, and the function topology will be discussed in Section~\ref{perm_section}.) Various authors subsequently have proved results of this form for other completely metrizable groups (e.g., see \cite{Bhattarcharjee} and \cite{GMR}). Then, Gartside and Knight~\cite{GK} undertook a more general study of the phenomenon in question. They defined a Polish group $G$ (i.e., one that is completely metrizable and separable) to be \emph{almost free} if the set $\{(g_1, \dots, g_n) \in G^n : \{g_1, \dots, g_n\} \ \text{freely generates a free subgroup of} \ G\}$ is comeagre in $G^n$ for each integer $n \geq 2$. They then gave various other characterizations of almost free groups and constructed new examples.

Our goal in this note is to extend the above definition of ``almost free" to groups that are not necessarily Polish (this is complicated by the fact that ``comeagre" is not a particularly useful notion for topological groups that are not completely metrizable, or more generally, ones that are not Baire spaces), generalize various results about almost free Polish groups to this context, and give new examples of groups with this property. Thus, given any infinite topological group $G$ and a cardinal $\kappa > 0$, we shall say that $G$ is \emph{almost $\kappa$-free} if the set $\{(g_i)_{i \in \kappa} \in G^\kappa : \{g_i\}_{i \in \kappa} \ \text{freely generates a free subgroup of} \ G\}$ is dense in $G^\kappa$. It turns out that if $G$ is a Polish group, then it is almost free (in the sense of Gartside and Knight) if and only if it is almost $n$-free for each positive integer $n$ (Lemma~\ref{def-equiv}). Therefore, our definition indeed generalizes the existing one, while preserving the intuition behind it. Namely, it still describes groups where sequences of elements that freely generate free subgroups are ubiquitous. 

We shall prove that an infinite topological group is almost free (i.e., almost $n$-free for each positive integer $n$) if and only if it is almost $\aleph_0$-free (Proposition~\ref{count-free}). We shall then show that if $G$ is a non-discrete Hausdorff topological group that contains a dense free subgroup of rank $\kappa > 0$, then $G$ is almost $\kappa$-free, as well as a more general version of this statement (Corollary~\ref{discrete-free2} and Theorem~\ref{discrete-free}). These are generalizations of parts of the main result of~\cite{GK}, though we employ very different methods (the proofs in~\cite{GK} are primarily topological, whereas ours are primarily algebraic).

Using the results above, we shall deduce that for any infinite set $\Omega$, the group $\Sym(\Omega)$ of all permutations of $\Omega$ is almost $2^{|\Omega|}$-free (with respect to the function topology), generalizing the aforementioned theorem of Dixon (Corollary~\ref{perm}). In a similar vein, we shall show that many $\kappa$-fold products of finite permutation groups are almost $2^\kappa$-free (Theorem~\ref{prod-main}). We shall also show that for any cardinal $\kappa > 0$, a free group of rank $\kappa$ is almost $\kappa$-free, with respect to any nondiscrete topology (Proposition~\ref{free-groups}). It turns out that every dense subgroup of a connected semi-simple real Lie group is also almost free (Corollary~\ref{lie}). Finally, we shall show that the direct product of two topological groups is almost $\kappa$-free if and only if at least one of the two groups is itself almost $\kappa$-free (Proposition~\ref{prod-example}). All the results mentioned in this paragraph give rise to examples of almost ($\kappa$-)free groups that are not Polish.

\subsection*{Acknowledgements}

The author is grateful to Yair Glasner for very helpful discussions about this material, to George Bergman for comments on an earlier draft of this note that have led to vast improvements, and to the referee for mentioning related literature.

\section{Generalities and examples} \label{generalities}

Throughout, $\Z_+$ will denote the set of positive integers. Given a set $\Gamma$ and a cardinal $\kappa$, we shall denote the direct product $\prod_\kappa \Gamma$ by $\Gamma^\kappa$. If $\Gamma$ is a topological space, then $\Gamma^\kappa$ will be understood to be a topological space under the product topology. Topological groups will not be assumed to be Hausdorff, unless it is specifically stated otherwise.

\begin{definition}
Given an infinite topological group $G$ and a cardinal $\kappa > 0$, let $G_\kappa$ denote the set $\, \{(g_i)_{i \in \kappa} \in G^\kappa : \{g_i\}_{i \in \kappa} \ \text{freely generates a free subgroup of} \ G\}$. We shall say that $G$ is \emph{almost $\kappa$-free} if $G_\kappa$ is dense in $G^\kappa$. Moreover, we shall say that $G$ is \emph{almost free} if $G$ is almost $n$-free for each $n \in\Z_+$.
\end{definition}

For the convenience of the reader, let us next recall some topological notions. A topological space is said to be \textit{Polish} if it is completely metrizable (i.e., it admits a metric with respect to which it is complete) and separable (i.e., it contains a countable dense subset). A subset of a topological space is called \textit{nowhere dense} if its closure contains no open subsets. Also, a subset is called \textit{comeagre} if it is the complement of a countable union of nowhere dense sets.

In~\cite{GK} Gartside and Knight defined a Polish topological group $G$ to be \emph{almost free} if for each $n \geq 2$, $G_n$ is comeagre in $G^n$, and they define $G$ to be \emph{almost countably free} if $G_{\aleph_0}$ is comeagre in $G^{\aleph_0}$. At first glance it might appear that our definition clashes with this one.  However, it turns out that if $G$ is a Polish (or even just a completely metrizable) group, then it is almost free in the sense of Gartside and Knight if and only if it is almost free in our sense, and $G$ is almost countably free if and only if it is almost $\aleph_0$-free, as the next lemma shows. Thus, from now on we shall use ``almost countably free" and ``almost $\aleph_0$-free" interchangeably.

\begin{definition}
Let $G$ be a group. Given a cardinal $\kappa > 0$, a set of variables $\, \{x_i\}_{i\in \kappa}$, and a free word $w = w(x_{i_1}, \dots, x_{i_n})$ $(i_1, \dots, i_n \in \kappa)$, let $C^\kappa (w)$ denote the set $\, \{(g_i)_{i \in \kappa} \in G^\kappa : w(g_{i_1}, \dots, g_{i_n}) = 1\}$.
\end{definition}

\begin{lemma} \label{def-equiv}
Let $G$ be a completely metrizable topological group, and let $\kappa$ be a cardinal satisfying $1 \leq \kappa \leq \aleph_0$. Then $G_\kappa$ is dense in $G^\kappa$ if and only if $G_\kappa$ is comeagre in $G^\kappa$.
\end{lemma}

\begin{proof}
Since $G$ is completely metrizable and $\kappa$ is a countable cardinal, $G^\kappa$ is completely metrizable as well. Hence, the ``if" direction follows from the Baire Category Theorem.

For the converse, we begin by noting that $G^\kappa \setminus G_\kappa = \bigcup_w C^\kappa(w)$. Since there are only countably many such words $w$, this is a union of countably many sets. Since the ``evaluation" map $G^n \rightarrow G$ induced by a word $w$ (of length $n$) is continuous, $\{1\}$ is closed in $G$ (as a Hausdorff space), and $C^\kappa(w) = w^{-1}(\{1\})$, we conclude that $C^\kappa(w)$ is closed in $G$. Since $G_\kappa$ is dense in $G^\kappa$, and since $C^\kappa(w) \subseteq G^\kappa \setminus G_\kappa$, it follows that $C^\kappa(w)$ cannot contain an open subset, for any free word $w$. Thus, $G^\kappa \setminus G_\kappa$ is the (countable) union of the nowhere dense sets $C^\kappa(w)$, and hence $G_\kappa$ is comeagre in $G^\kappa$.
\end{proof}

Let us pause to give some examples. As mentioned in the Introduction, Dixon~\cite{Dixon} showed that $\Sym(\Z_+)$ is almost free (see also~\cite{GMR}; more on this below). In~\cite{GK} Gartside and Knight gave a number of other examples of Polish almost free groups. Specifically, all Polish oligomorphic permutation groups have this property (i.e., subgroups of the group of all permutations of an infinite set $\Omega$, $\Sym(\Omega)$, whose action on $\Omega^n$ has only finitely many orbits, for all $n \in \Z_+$). Every finite-dimensional connected non-solvable Lie group is almost free (see also~\cite{Epstein}), as is the absolute Galois group of the rational numbers (i.e., the group of all automorphisms of the algebraic closure of $\Q$). Although such a group is completely metrizable rather than Polish, an inverse limit of wreath products of nontrivial groups is also almost free (see also~\cite{Bhattarcharjee}).

The next result allows us to construct new almost free groups that are not necessarily completely metrizable, from existing ones.

\begin{proposition} \label{prod-example}
Let $\kappa > 0$ be a cardinal, and let $G$ and $H$ be topological groups. Then the direct product $G \times H$ $($regarded as a topological group in the product topology$)$ is almost $\kappa$-free if and only if at least one of $G$ and $H$ is almost $\kappa$-free.
\end{proposition}

\begin{proof}
Suppose that $G$ is almost $\kappa$-free, and let $U \subseteq (G \times H)^\kappa$ be a nonempty open set. We wish to show that $U \cap (G \times H)_\kappa \neq \emptyset$. After passing to a subset and reindexing, we may assume without loss of generality that $U$ is of the form $$(U_1 \times V_1) \times \ldots \times (U_n \times V_n) \times \prod_{\kappa \setminus \{1, \ldots, n\}} (G \times H),$$ where $U_1, \ldots, U_n \subseteq G$, $V_1, \ldots, V_n \subseteq H$ are some nonempty open sets. Since $G$ is almost $\kappa$-free, we can find a sequence $(g_i)_{i \in \kappa} \in G_\kappa$ such that $g_i \in U_i$ for $i \in \{1, \ldots, n\}$. Let $(h_i)_{i \in \kappa} \in H^\kappa$ be any sequence such that $h_i \in V_i$ for $i \in \{1, \ldots, n\}$. Setting $f_i = (g_i, h_i) \in G\times H$ for all $i\in \kappa$, we have $(f_i)_{i \in \kappa} \in U$. Further, since $\{g_i\}_{i \in \kappa}$ freely generates a free subgroup, so does $\{f_i\}_{i \in \kappa}$. (If the $f_i$ satisfy some free word, then so do their projections onto $G$.) Hence $(f_i)_{i \in \kappa} \in U \cap (G \times H)_\kappa$, as desired. Similarly, if $H$ is almost $\kappa$-free, then so is $G \times H$.

For the converse, suppose that neither $G$ nor $H$ is almost $\kappa$-free. Then we can find open subsets $U = \prod_{i \in \kappa} U_i \subseteq G^\kappa$ and $V = \prod_{i \in \kappa} V_i \subseteq H^\kappa$ such that $U \cap G_\kappa = \emptyset = V \cap H_\kappa$. Set $W = \prod_{i \in \kappa} (U_i \times V_i) \subseteq (G \times H)^\kappa,$ and let $(f_i)_{i \in \kappa} = ((g_i, h_i))_{i \in \kappa} \in W$ be an arbitrary element. We wish to show that  $(f_i)_{i \in \kappa} \notin (G \times H)_\kappa$, which implies that $G \times H$ is not almost $\kappa$-free, since $W$ is a nonempty open subset of $(G \times H)^\kappa$. After reindexing, if necessary, we can find free words $w_G = w_G(x_1, \dots, x_n)$ and $w_H = w_H(x_1, \dots, x_n)$, for some $n \geq 1$, such that $w_G(g_1, \dots, g_n) = 1$ and $w_H(h_1, \dots, h_n) = 1$, by our assumptions on $U$ and $V$. It follows that $w_Gw_Hw_G^{-1}w_H^{-1}(f_1, \dots, f_n) = 1$. If $w_Gw_Hw_G^{-1}w_H^{-1}$ is not the trivial word, then this shows that $(f_i)_{i \in \kappa} \notin (G \times H)_\kappa$, as desired. If, however, $w_Gw_Hw_G^{-1}w_H^{-1}$ is the trivial word, then it must be the case that $w_G$ and $w_H$ are powers of the same element. (This follows from the fact that all subgroups of a free group are free, and hence every commutative subgroup is cyclic.) If $n = 1$, then it must be the case that there is a (nontrivial) free word of the form $x^m$ ($m \in \Z_+$) satisfied by $f_1$. Otherwise, assuming that $n > 1$, we can find a free word $w = w(x_1, \dots, x_n)$ that commutes with neither $w_G$ nor $w_H$. Also, we still have $ww_Hw^{-1}w_H^{-1}(h_1, \dots, h_n) = 1$. Therefore, replacing $w_H$ with $ww_Hw^{-1}w_H^{-1}$ in the expression $w_Gw_Hw_G^{-1}w_H^{-1}$, we obtain a nontrivial free word satisfied by $f_1, \dots, f_n$.
\end{proof}

Here is another way of constructing new almost free groups from old ones, which will be needed in the sequel.

\begin{lemma} \label{dense-subgroup}
Let $\kappa > 0$ be a cardinal, and let $G$ be a topological group containing a dense subgroup $H$ which is almost $\kappa$-free with respect to the induced topology. Then $G$ is itself almost $\kappa$-free. 
\end{lemma}

\begin{proof}
Let $U \subseteq G^\kappa$ be a nonempty open subset. Then $U \cap H^\kappa$ is a nonempty open subset of $H^\kappa$, and hence $\emptyset  \neq (U \cap H^\kappa) \cap H_\kappa = U \cap H_\kappa \subseteq U \cap G_\kappa$. Therefore, $G$ is almost $\kappa$-free.
\end{proof}

Let us next record another useful basic fact.

\begin{lemma} \label{down}
Let $G$ be a topological group, and let $\, 0 < \lambda < \kappa$ be cardinals. If $G$ is almost $\kappa$-free, then it is almost $\lambda$-free.
\end{lemma}

\begin{proof}
Writing $\kappa = \lambda + \nu$, for some cardinal $\nu$, we may identify $G^\kappa$ with $G^\lambda \times G^\nu$. Let $U \subseteq G^\lambda$ be a nonempty open subset. Then $V = U \times G^\nu \subseteq G^\kappa$ is an open subset as well. Since $G$ is almost $\kappa$-free, $G_\kappa \cap V \neq \emptyset$. Let $\pi_\lambda : G^\kappa \rightarrow G^\lambda$ be the natural projection with respect to the product decomposition $G^\kappa = G^\lambda \times G^\nu$. Then $\pi_\lambda(G_\kappa) \cap U = \pi_\lambda(G_\kappa) \cap \pi_\lambda (V) \neq \emptyset$. Now, $\pi_\lambda(G_\kappa) \subseteq G_\lambda$, since if a sequence of elements freely generates a free group, then so does any subsequence. Hence, $G_\lambda \cap U \neq \emptyset$, showing that $G$ is almost $\lambda$-free.
\end{proof}

The following is a generalization of a result of Gartside and Knight in~\cite{GK} about Polish groups.

\begin{proposition} \label{count-free}
Let $G$ be a topological group. Then $G$ is almost free if and only if $G$ is almost countably free.
\end{proposition}

\begin{proof}
The ``if" direction follows immediately from the previous lemma. For the converse, let $U \subseteq G^{\aleph_0}$ be a nonempty open set. After passing to a subset and reindexing, we may assume without loss of generality that $U$ is of the form $U_1 \times \ldots \times U_n \times G \times G \times \ldots$, where $U_1, \ldots, U_n \subseteq G$ are some nonempty open sets. Since $G$ is almost free, we can find some $(g_1, \ldots, g_{n+2}) \in U_1 \times \ldots \times U_n \times G \times G$ such that $g_1, \ldots, g_{n+2}$ freely generate a free subgroup of $G$. Then, in particular, $g_{n+1}, g_{n+2}$ also freely generate a free subgroup of $G$. It is well known that one can embed a free group on $\aleph_0$ generators into a free group on two generators. Thus, we can find $\{h_1, h_2, \ldots\} \subseteq \langle g_{n+1}, g_{n+2} \rangle$ (the subgroup generated by $g_{n+1}$ and $g_{n+2}$) that freely generates a free subgroup. Then $\{g_1, \ldots, g_n\} \cup \{h_1, h_2, \ldots \}$ also freely generates a free subgroup of $G$. (If there is a free word $w(x_1, \dots, x_m)$ and elements $k_1, \ldots, k_m \in \{g_1, \ldots, g_n\} \cup \{h_1, h_2, \ldots \}$ such that $w(k_1, \dots, k_m) = 1$, then there must be a free word $\bar{w}(x_1, \dots, x_{n+2})$ such that $\bar{w}(g_1, \dots, g_{n+2}) = 1$.) Thus $(g_1, \ldots, g_n, h_1, h_2, \ldots) \in G_{\aleph_0} \cap U$, showing that $G$ is almost countably free.
\end{proof}

When studying almost $\kappa$-freeness of groups, it is harmless to focus on ones that are not discrete, as the next observation shows.

\begin{lemma} \label{discrete}
Let $G$ be a discrete topological group. Then $G$ is not almost $\kappa$-free for any cardinal $\kappa > 0$.
\end{lemma}

\begin{proof}
It suffices to note that $\{1\} \times \prod_{\kappa \setminus \{1\}} G \subseteq G^\kappa \setminus G_\kappa$ is an open set, for any cardinal $\kappa > 0$.
\end{proof}

\section{Dense free subgroups}

The main goal of this section is to prove that, under mild assumptions, a group containing a dense free subgroup of rank $\kappa > 0$ is almost $\kappa$-free (Theorem~\ref{discrete-free} and Corollary~\ref{discrete-free2}). We begin with several lemmas.

\begin{lemma} \label{non_comm}
Let $F$ be a non-discrete free group of rank $\geq 2$, and let $U \subseteq F$ be an open neighborhood of the identity. Then $U$ is not commutative. Hence, one can find two elements of $U$ that freely generate a free subgroup.
\end{lemma}

\begin{proof}
Suppose that $U$ is commutative, and let $S \subseteq F$ be the subgroup of all elements that commute with all elements of $U$. Since $F$ is free and non-discrete, there is an element $x \in F \setminus \{1\}$ such that $U \subseteq S = \langle x \rangle$. (As in the proof of Proposition~\ref{prod-example}, this follows from the fact that all subgroups of $F$ are free, and hence every commutative subgroup is cyclic.) Since $F$ has rank $\geq 2$, we can find an element $y \in F$ such that $x$ and $y$ freely generate a free subgroup. Now, since $x \cdot 1 \cdot x^{-1}, y \cdot 1 \cdot y^{-1} \in U$, we can find a neighborhood of the identity $W$ such that $xWx^{-1}, yWy^{-1} \subseteq U$. It follows that $W \subseteq \langle x \rangle$. (For, if $g \in W$ does not commute with an element of $U$, then neither does $xgx^{-1} \in U$.) Since $F$ is non-discrete, $W \neq \{1\}$, and hence $yWy^{-1}$ must contain an element that does not commute with $x$ (since $x$ and $y$ freely generate a free subgroup), contradicting the assumption that $U$ is commutative. Therefore, $U$ cannot be commutative.

The final claim follows from the aforementioned fact that all subgroups of $F$ are free, and hence if two elements of $F$ generate a noncommutative subgroup, then it must have rank $2$.
\end{proof}

\begin{lemma} \label{prod-neigh}
Let $F$ be a topological group, let $U_1, \ldots, U_n \subseteq F$ be nonempty open subsets, let $g_1 \in U_1, \ldots, g_n \in U_n$ be any elements, and let $m \in \Z_+$. Then there is an open neighborhood $U$ of the identity such that for all $\, 0 \leq j \leq m$ and $\, 1 \leq i \leq n$, we have $U^{m-j}g_i U^j \subseteq U_i$.
\end{lemma}

\begin{proof}
For each $0 \leq j \leq m$, let $w_j = w_j(x,y)$ be the free word $x^{m-j}yx^j$. Since the ``evaluation" map $F^{m+1} \rightarrow F$ induced by the word $w_j$ is continuous, and since $1^{m-j}g_i1^j \in U_i$, for each $0 \leq j \leq m$ and $1 \leq i \leq n$ there is an open neighborhood $U_{ij}$ of the identity such that $U_{ij}^{m-j}g_i U_{ij}^j \subseteq U_i$. Then the open neighborhood of the identity $U = \bigcap_{i,j} U_{ij}$ has the desired properties.
\end{proof}

\begin{lemma} \label{fin-case}
Let $F$ be a non-discrete free group of rank $\kappa \geq 2$. Then $F$ is almost free.
\end{lemma}

\begin{proof}
Let $\{x_i\}_{i \in \kappa}$ be a set of free generators for $F$. Given an element $g \in F$, let $L(g)$ denote the length of $g$ as a reduced monoid word in $\{x_i\}_{i \in \kappa} \cup \{x_i^{-1}\}_{i \in \kappa}$. Also, let $n$ be a positive integer, and let $U_1, \ldots, U_n \subseteq F$ be nonempty open subsets. We shall construct elements $f_1 \in U_1, \ldots, f_n \in U_n$ that freely generate a free subgroup.

We begin by picking arbitrary elements $g_1 \in U_1, \ldots, g_n \in U_n$. Let $m$ be an integer greater than $2L(g_i) + 2n + 2$, for all $1 \leq i \leq n$, and let $U$ be a neighborhood of the identity as in Lemma~\ref{prod-neigh}. By Lemma~\ref{non_comm}, there are elements $y,z \in U$ that freely generate a free subgroup. For each $1 \leq i \leq n$, we can find elements $h_{i,1}, h_{i,2} \in \langle y, z \rangle$ such that $f_i = z^{-i}yh_{i,1}g_ih_{i,2}yz^i$ when reduced, yields a word of the form $z^{-i}ywyz^i$, where $w$ is reduced as a word in the $x_i$, and where no additional reductions are possible other than those resulting from ones that occur in $z^{-i}y$ or $yz^i$ (as words in the $x_i$). Moreover, we may choose the $h_{i,1}$ and $h_{i,2}$ so that each has length no greater than $L(g_i)$ as a monoid word in $\{y, y^{-1}, z, z^{-1}\}$. (For, $g_i$ cannot cancel more than $L(g_i)$ copies of $y$ and $z$ on either side.) Hence, by our choice of $U$, for each $1 \leq i \leq n$, we have $f_i \in U_i$. 

Now, let $z^{i_1}y^{i_2}ry^{i_3}z^{i_4}$ and $z^{i_5}y^{i_6}sy^{i_7}z^{i_8}$ be two reduced group words in the $x_i$ (except for possible unresolved reductions in sub-words of the form $z^jy^k$ at the beginning and sub-words of the form $y^kz^j$ at the end), for some $r,s \in F$ and nonzero integers $i_1, \ldots, i_8$, where $i_4 \neq - i_5$. Then $$(z^{i_1}y^{i_2}ry^{i_3}z^{i_4})(z^{i_5}y^{i_6}sy^{i_7}z^{i_8}) = z^{i_1}y^{i_2}(ry^{i_3}z^{i_4+i_5}y^{i_6}s)y^{i_7}z^{i_8},$$ when reduced, also yields a reduced group word (with the same caveat) of the same form. It follows that if $w = w(t_1, \ldots, t_n)$ is a reduced free group word, then $w(f_1, \ldots, f_n) \neq 1$, and hence $f_1, \ldots, f_n$ freely generate a free subgroup.
\end{proof}

We are now in a position to characterize when a free group is almost $\kappa$-free.

\begin{proposition} \label{free-groups}
Let $F$ be a free topological group of rank $\kappa$, for some cardinal $\kappa > 0$. Then the following are equivalent.
\begin{enumerate}
\item[$(1)$] The topology on $F$ is non-discrete.
\item[$(2)$] $F$ is almost $\kappa$-free.
\end{enumerate}
Moreover, in the above situation, if $\, 2 \leq \kappa \leq \aleph_0$, then $F$ is almost countably free.
\end{proposition}

\begin{proof}
By Lemma~\ref{discrete}, (2) implies (1). Let us show that (1) implies (2) and the final claim.

First, suppose that $\kappa = 1$. Since $F$ is non-discrete, for any nonempty open subset $U \subseteq F$, we can find an element $f \in U\setminus \{1\}$, and $f$ freely generates a free subgroup. Hence $F$ is almost $1$-free. 

Next, suppose that $2 \leq \kappa \leq \aleph_0$. By Lemma~\ref{fin-case}, $F$ is almost free, and therefore, Proposition~\ref{count-free} implies that $F$ is almost countably free. In particular, $F$ is almost $\kappa$-free, by Lemma~\ref{down}.

Finally, suppose that $\aleph_0 \leq \kappa$, and let $\{x_i\}_{i \in \kappa}$ be a set of free generators for $F$. Also, let $U \subseteq F^\kappa$ be a nonempty open subset. We wish to show that $U \cap F_\kappa \neq \emptyset$. After passing to a subset and reindexing, we may assume without loss of generality that $U$ is of the form $U_1 \times \ldots \times U_n \times \prod_{\kappa \setminus \{1, \ldots, n\}} F$, where $U_1, \ldots, U_n \subseteq F$ are some nonempty open sets. By Lemma~\ref{fin-case}, we can find elements $f_1 \in U_1, \ldots, f_n \in U_n$ that freely generate a free subgroup of $F$. As words in $\{x_i\}_{i \in \kappa}$, the elements $f_1, \ldots,  f_n$ involve only finitely many of the $x_i$, say, $\{x_{i_1}, \ldots, x_{i_m}\}$ $(i_1, \ldots, i_m \in \kappa)$. Let $\Lambda = \kappa \setminus \{x_{i_1}, \ldots, x_{i_m}\}$. Then $\{f_1, \ldots, f_n\} \cup \{x_i\}_{i \in \Lambda}$ freely generates a free subgroup of rank $\kappa$. Writing $\{f_i\}_{i \in \kappa, \, n<i} = \{x_i\}_{i \in \Lambda}$, we therefore have $(f_i)_{i\in \kappa} \in U \cap F_\kappa$, as desired.
\end{proof}

We recall that no countable non-discrete group $G$ can be completely metrizable. For, unless there is an element $g \in G$ that is isolated, $G$ can be written as a countable union of nowhere dense sets, namely the singleton sets, which contradicts the Baire Category Theorem. On the other hand, if some $g \in G$ is isolated, then $G$ must be discrete. Hence, in particular, Proposition~\ref{free-groups} gives additional examples of non-Polish almost $\kappa$-free groups when $\kappa$ is countable.

We are ready for our main result.

\begin{theorem} \label{discrete-free}
Let $G$ be a topological group, and let $\kappa > 0$ be a cardinal. Suppose that $G$ contains a free subgroup $F$ of rank $\kappa$, such that $F\setminus \{1\}$ is dense. Then $G$ is almost $\kappa$-free. Moreover, if $\, 2 \leq \kappa \leq \aleph_0$, then $G$ is almost countably free.
\end{theorem}

\begin{proof}
Our hypotheses imply that the induced topology on $F$ is non-discrete. The result now follows from Proposition~\ref{free-groups} and Lemma~\ref{dense-subgroup}.
\end{proof}

The following consequence of Theorem~\ref{discrete-free} generalizes a result of Gartside and Knight in~\cite{GK} about Polish groups.

\begin{corollary} \label{discrete-free2}
Let $G$ be a non-discrete Hausdorff topological group, and let $\kappa > 0$ be a cardinal. Suppose that $G$ contains a dense free subgroup $F$ of rank $\kappa$. Then $G$ is almost $\kappa$-free. Moreover, if $\, 2 \leq \kappa \leq \aleph_0$, then $G$ is almost countably free.
\end{corollary}

\begin{proof}
Since $G$ is non-discrete, it contains a point that is not isolated. By translation-invariance of the topology, it follows that $1$ is not isolated. Since $G$ is Hausdorff, this implies that every nonempty open set contains a nonempty open subset which is not a neighborhood of the identity. Thus, if $F \subseteq G$ is dense, then $F \setminus \{1\}$ must be dense as well. The result now follows from Theorem~\ref{discrete-free}.
\end{proof}

By Lemma~\ref{discrete}, the non-discreteness assumption in the above corollary is necessary. Using similar reasoning, it is possible to construct an example showing that the Hausdorff condition is necessary as well. For, let $F$ be a discrete free group of rank $\kappa > 0$, and let $H$ be an indiscrete group which contains no nontrivial free subgroups (e.g., a group where all elements have finite order). Then $F \times \{1\}$ is a dense free subgroup of $G = F \times H$, which is clearly a non-discrete non-Hausdorff group (if $H \neq \{1\}$). However, $G$ cannot be almost $\kappa$-free, since $(\{1\} \times H) \times \prod_{\kappa \setminus \{1\}} (F \times H)$ is an open subset of $G^\kappa \setminus G_\kappa$.

Garside and Knight~\cite{GK} showed that if $G$ is a non-discrete Polish topological group that is almost countably free, then it contains a dense free subgroup (of rank $\aleph_0$). Thus, one might wonder whether more general converses to Theorem~\ref{discrete-free} and Corollary~\ref{discrete-free2} hold. Let us note that it is possible to construct groups that are almost $\kappa$-free but have no dense free subgroups at all. For instance, let $\kappa < \lambda$ be infinite cardinals, let $F$ be a free group of rank $\kappa$ (with any non-discrete topology; e.g., the profinite topology), and let $A$ be a discrete abelian group of cardinality $\lambda$. Then, by Proposition~\ref{prod-example}, $F \times A$ is almost $\kappa$-free. However, since $A$ is abelian, all free subgroups of $F \times A$ must have cardinality $\leq \kappa$, and hence cannot be dense.

One can even construct a similar example of a non-discrete Hausdorff group that is almost $\kappa$-free and has a dense subgroup on $\kappa$ generators but no dense free subgroups. For instance, let $F$ be a non-discrete free group of rank $\aleph_0$, and let $G = \bigoplus_{\aleph_0} F$, the subgroup of $\prod_{\aleph_0} F$ consisting of elements having only finitely many nonidentity components. We regard $G$ as a topological group under the topology induced by the product topology on $\prod_{\aleph_0} F$ (which makes $G$ Hausdorff if $F$ is). Again, by Proposition~\ref{prod-example}, $G$ is almost countably free, and as a countable group it clearly has a dense subgroup on $\aleph_0$ generators. However, for each $j \in \aleph_0$, any dense subgroup of $G$ must contain an element of the form $g_j = (f_i)_{i \in \aleph_0}$, where $f_i = 1$ for all $i \neq j$, and $f_j \neq 1$. Now, any two such elements $g_j$ and $g_k$ commute, and they cannot be powers of the same element of $G$ unless $j = k$ (since only $1 \in F$ has finite order). Hence, no dense subgroup of $G$ can be free. Thus, the most obvious attempts to relax the complete metrizability assumption in the aforementioned converse to Corollary~\ref{discrete-free2} fail.

In~\cite{BG} Breuillard and Gelander showed that every dense subgroup of a connected semi-simple real Lie group $G$ contains a free subgroup of rank $2$ that is dense (in $G$). Hence, we have the following consequence of Corollary~\ref{discrete-free2}.

\begin{corollary} \label{lie}
Every dense subgroup of a connected semi-simple real Lie group is almost countably free.
\end{corollary}

Many such groups are not Polish; for instance, the countable ones.

We note in passing that the main results of this section have a flavor similar to that of the following theorem from~\cite{Gelander}.

\begin{theorem}[Gelander]
Let $G$ be a connected compact Lie group, let $n \geq 3$ be an integer, and let $H \subseteq G$ be an $(n-1)$-generated dense subgroup. Then $\, \{(h_1, \dots, h_n) \in H^n : \langle h_1, \dots, h_n \rangle = H \}$ is dense in $G^n$.
\end{theorem}

\section{Permutations} \label{perm_section}

Let $\Omega$ be an infinite set. Regarding $\Omega$ as a discrete topological space, the monoid $\Self(\Omega)$ of all set maps from $\Omega$ to itself becomes a topological space under the \emph{function topology}. A subbasis of open sets in this topology is given by the sets $\{f \in \Self(\Omega) : f(\alpha) = \beta\}$ $(\alpha,\beta \in \Omega).$ It is easy to see that composition of maps is continuous in this topology. The group $\Sym(\Omega)$ of all permutations of $\Omega$ inherits from $\Self(\Omega)$ the function topology. Moreover, when restricted to $\Sym(\Omega)$, this topology makes $( )^{-1}$ continuous, since $$\{f \in \Sym(\Omega) : f(\alpha) = \beta\}^{-1} = \{f \in \Sym(\Omega) : f(\beta) = \alpha\},$$ turning $\Sym(\Omega)$ into a topological group, which can easily be seen to be Hausdorff.

In the case where $\Omega = \Z_+$, one can put a metric $d$ on $\Sym(\Omega)$ which induces the function topology. Specifically, given $f, g \in \Sym(\Z_+)$, let $d(f,g) = 0$ if $f = g$, and otherwise let $d(f,g) = 2^{-n}$, where $n \in \Z_+$ is the least number such that either $f(n) \neq g(n)$ or $f^{-1}(n) \neq g^{-1}(n)$. A subbasis of open sets in the topology induced by this metric on $\Sym(\Z_+)$ consists of sets of the form 
$$(\ast) \ \ \bigcap_{i = 1}^n \{f \in \Sym(\Z_+) : f(i)  = \alpha_i \ \text{and} \ f^{-1}(i) = \beta_i\} \ (n, \alpha_1, \ldots, \alpha_n, \beta_1, \ldots, \beta_n \in \Z_+).$$ We note that, with $\alpha_i$ and $\beta_i$ as above, both $\{f \in \Sym(\Z_+) : f(i)  = \alpha_i\}$ and $\{f \in \Sym(\Z_+) : i = f(\beta_i)\}$ are open sets with respect to the function topology, as defined in the previous paragraph, and hence so is $$\{f \in \Sym(\Z_+) : f(i)  = \alpha_i \ \text{and} \ f^{-1}(i) = \beta_i\}$$ $$= \{f \in \Sym(\Z_+) : f(i)  = \alpha_i\} \cap \{f \in \Sym(\Z_+) : i = f(\beta_i)\}.$$ Thus, a set open in the topology induced by $d$ is open in the function topology. Conversely, a set of the form $\{f \in \Sym(\Z_+) : f(n) = m\}$ can be expressed as the union of all the sets of the form ($\ast$) where $\alpha_n = m$. Therefore, a set that is open in the function topology is open in the topology induced by $d$ as well, showing that the two are in fact the same topology.

We wish to show that for any infinite set $\Omega$, $\Sym(\Omega)$ is almost $2^{|\Omega|}$-free, with respect to the topology described above. This task can be accomplished easily if we rely on the following result from~\cite{Bruijn}.

\begin{theorem}[de Bruijn] \label{debruijn}
Let $\Omega$ be an infinite set. Then $\, \Sym(\Omega)$ contains a free subgroup of rank $\, 2^{|\Omega|}.$ Moreover, this free subgroup can be taken to be dense $($in the function topology$)$.
\end{theorem}

The last claim in the above theorem does not actually appear in~\cite{Bruijn} but does, however, follow from de Bruijn's proof, as noted by Hodges. See~\cite{MS} for a discussion of this, along with a model-theoretic generalization of the result.

Applying Corollary~\ref{discrete-free2} to Theorem~\ref{debruijn}, we obtain the following.

\begin{corollary} \label{perm}
Let $\Omega$ be an infinite set. Then $\, \Sym(\Omega)$ is almost $\, 2^{|\Omega|}$-free, with respect to the function topology. 
\end{corollary}

Corollary~\ref{perm} is a generalization of the result of Dixon~\cite{Dixon} mentioned in the Introduction, showing that $\Sym(\Z_+)$ is almost free with respect to the function topology. (See Lemma~\ref{down}. See also~\cite{BR} for generalizations of Dixon's result and Theorem~\ref{debruijn} in a different direction--to automorphisms groups of relatively free algebras.) We note that $\Sym(\Z_+)$ is complete with respect to the metric defined above, and it is not hard to see that it is separable as well. Thus, with respect to the function topology, $\Sym(\Omega)$ is Polish whenever $|\Omega| = \aleph_0$. This is not the case, however, if $|\Omega| > \aleph_0$. For, in this situation, $\Sym(\Omega)$ is neither separable nor metrizable, since it is not first-countable. (A topological space is said to be \emph{first-countable} if each point has a countable base for its system of neighborhoods. Every metric space is first-countable, since given a point $p$, the open balls centered at $p$ of radii $1/n$, for $n \in \Z_+$, form such a countable base for this point.)

As mentioned above, in~\cite{MS} Melles and Shelah proved a more general version of Theorem~\ref{debruijn}, showing that for certain models $M,$ the automorphism group $\Aut(M)$ of $M$ contains a free subgroup of rank $2^{|M|}$ that is dense (in the function topology). Hence, Corollary~\ref{perm} can be generalized accordingly. 

Let us next construct a different sort of almost $\kappa$-free group of permutations. For each $m \in \Z_+$, let $S_m = \Sym(\{1, 2, \ldots, m\})$ denote the group of all permutations of the set $\{1, 2, \ldots, m\}$. Let $\kappa$ be an infinite cardinal, let $\varphi : \kappa \rightarrow \Z_+$ be a function, and let $G = \prod_{i \in \kappa} S_{\varphi(i)}$. Endowing each $S_{\varphi(i)}$ with the discrete topology, $G$ becomes a Hausdorff topological group in the product topology. As with $\Sym(\Omega)$ above, $G$ cannot be Polish if $\kappa > \aleph_0$. We shall show that $G$ is almost $2^\kappa$-free with respect to this topology, assuming that for each $l \in \Z_+$, $|\{i : \varphi(i) \geq l\}| = \kappa$. The argument is divided into several steps.

\begin{lemma} \label{prod1}
Let $w = w(x_1, \ldots, x_n)$ be a free word. Then there exist $m \in \Z_+$ and permutations $f_1, \ldots, f_n \in S_m$ such that $w(f_1, \ldots, f_n) \neq 1$.
\end{lemma}

\begin{proof}
Since $\Sym(\Z_+)$ contains a free group of rank $2^{\aleph_0}$, we can find $g_1, \ldots, g_n \in \Sym(\Z_+)$ such that $w(g_1, \ldots, g_n) \neq 1$. (Actually, here we only need the fact that $\Sym(\Z_+)$ contains a free group of rank $n$, and this is easy to see without reference to Theorem~\ref{debruijn}. For, every countable group acts on itself and hence can be embedded in $\Sym(\Z_+)$. This, in particular, holds for countable free groups.) Thus, for some distinct $i, j \in \Z_+$ we have $w(g_1, \ldots, g_n)(i) = j$. Write $w(g_1, \ldots, g_n) = h_1 \ldots h_k$, where $h_1, \ldots, h_k \in \{g_1, \ldots, g_n\} \cup \{g_1^{-1}, \ldots, g_n^{-1}\}$. Let $\Gamma = \{i, h_k(i), h_{k-1}h_k(i), \ldots, h_1 \dots h_k(i)\}$, and set $$\Delta = \Gamma \cup g_1(\Gamma) \cup \dots \cup g_n(\Gamma) \cup g_1^{-1}(\Gamma) \cup \dots \cup g_n^{-1}(\Gamma).$$
Let $m$ be the maximal element of the finite set $\Delta \subseteq \Z_+$. Then we can find $f_1, \ldots, f_n \in S_m$ such that for each $i \in \{1, \ldots, n\}$, $g_i$ agrees with $f_i$ on $\Gamma$, and $g_i^{-1}$ agrees with $f_i^{-1}$ on $\Gamma$. In particular, we have $w(f_1, \ldots, f_n)(i) = j \neq i$, as desired.
\end{proof}

\begin{proposition} \label{prod2}
Let $\kappa$ be an infinite cardinal, let $\varphi : \kappa \rightarrow \Z_+$ be a function, and let $G = \prod_{i \in \kappa} S_{\varphi(i)}$. Suppose that $\, |\{i : \varphi(i) \geq l\}| = \kappa$ for each $l \in \Z_+$. Then $G$ contains a free subgroup of rank $\kappa$.
\end{proposition}

\begin{proof}
First, suppose that $\kappa = \aleph_0$, and let us show that for all $n \in \Z_+$, $G$ contains a free subgroup of rank $n$. Let $\{w_i\}_{i \in \aleph_0}$ be an enumeration of all the free words in the $n$ letters $x_1, \ldots, x_n$. By Lemma~\ref{prod1}, to each $i \in \aleph_0$ we can assign an element $m_i \in \aleph_0$ such that there exist permutations $f_{i1}, \ldots, f_{in} \in S_{\varphi(m_i)}$ for which $w_i(f_{i1}, \ldots, f_{in}) \neq 1$. Moreover, by our hypotheses on $G$, we may assume that $m_i \neq m_j$ for $i \neq j$. Now, let $g_1, \ldots, g_n \in G$ be any permutations such that the natural projection of $g_j$ on $S_{\varphi(m_i)}$ is $f_{ij}$ ($1 \leq j \leq n$, $i \in \aleph_0$). Then for all $i \in \aleph_0$ we have $w_i(g_1, \ldots, g_n) \neq 1$, since the natural projection of $w_i(g_1, \ldots, g_n)$ on $S_{\varphi(m_i)}$ is $w_i(f_{i1}, \ldots, f_{in})$. Hence, $g_1, \ldots, g_n$ freely generate a free group.

For the general case, upon relabeling (and using the axiom of choice), we can identify $G$ with $\prod_{i \in \kappa} H_i$, where $H_i = \prod_{j \in \aleph_0} S_{\varphi_i(j)}$, each $\varphi_i : \aleph_0 \rightarrow \Z_+$ ($i \in \kappa$) is a function, and we have $|\{j : \varphi_i(j) \geq l\}| = \aleph_0$ for all $l \in \Z_+$ and $i \in \kappa$. We shall construct a set $\{f_i\}_{i \in \kappa} \subseteq G$ that freely generates a free subgroup. Let $\{U_i\}_{i \in \kappa}$ be a well-ordering of all the finite subsets of $\{f_i\}_{i \in \kappa}$ (for now, treated simply as a set of symbols), and write $U_i = \{f_{i1}, \ldots, f_{in_i}\}$, where $n_i = |U_i|$. Now, we define the permutations $f_i$ so that, for each $i \in  \kappa$, the natural projections of $f_{i1}, \ldots, f_{in_i}$ on $H_i$ generate a free subgroup (this is possible, by the previous paragraph), and if $f_i \notin U_j$ for some $i, j \in \kappa$, then we let the natural projection of $f_i$ on $H_j$ be the identity. It follows that if $\{g_1, \ldots, g_n\} \subseteq \{f_i\}_{i \in \kappa}$ is any finite subset, then $g_1, \ldots, g_n$ cannot satisfy a free word, and hence $\{f_i\}_{i \in \kappa}$ freely generates a free group.
\end{proof}

The above proposition actually implies that the group $G$ in question contains a free subgroup of rank $2^\kappa$. To show this we shall need the following result from~\cite{Bergman}. (The second paragraph of the proof of Proposition~\ref{prod2} employs a similar idea to that used to prove this theorem.)

\begin{theorem}[Bergman] \label{bergman}
Let $\, V$ be any variety of finitary $($general$)$ algebras, let $F$ be a free $\, V$-algebra on $\, \aleph_0$ generators, and let $\kappa$ be an infinite cardinal. Then $F^\kappa$ contains a free $\, V$-algebra on $2^\kappa$ generators.

In particular, if $F$ is a free group on $\, \aleph_0$ generators, then $F^\kappa$ contains a free group on $2^\kappa$ generators.
\end{theorem}

The proof of the following corollary is similar to that of Theorem~\ref{debruijn}.

\begin{corollary} \label{prod3}
Let $\kappa$ and $G$ be as in Proposition~\ref{prod2}. Then $G$ contains a free subgroup of rank $2^\kappa$.
\end{corollary}

\begin{proof}
First, we note that $G$ contains a $\kappa$-fold direct product of groups of the same form, by the same argument as in the second paragraph of the proof of Proposition~\ref{prod2}. (For, $\kappa$ can be written as a disjoint union of $\kappa$ subsets of cardinality $\kappa$.) Now, by the same proposition, it follows that $G$ contains a $\kappa$-fold direct product of free groups of rank $\aleph_0$. Thus, by Theorem~\ref{bergman}, $G$ contains a free group of rank $2^\kappa$.
\end{proof}

We require one more ingredient to construct our desired almost $\kappa$-free group of permutations.

\begin{lemma} \label{prod4}
Let $\kappa$ be an infinite cardinal, and let $G = \prod_{i \in \kappa} G_i$, where each $G_i$ is of the form $S_m$, for some $m \in \Z_+$. Suppose that $f_1, \ldots, f_n \in G$ freely generate a free subgroup, and let $g_1, \ldots, g_n \in G$ be any permutations such that there exists a finite subset $\, \Gamma \subseteq \kappa$ with the property that for each $j \in \{1, \ldots, n\}$ and $i \in \kappa \setminus \Gamma$, the natural projections of $g_j$ and $f_j$ on $G_i$ are equal. Then $g_1, \ldots, g_n$ also freely generate a free subgroup.
\end{lemma}

\begin{proof}
Suppose, on the contrary, that $w(g_1, \ldots, g_n) = 1$ for some free word $w$. Then the natural projection of $w(f_1, \ldots, f_n)$ on $\prod_{i \in \kappa \setminus \Gamma} G_i$ is also the identity. Now, $\prod_{i \in \Gamma} G_i$ is a finite group, and hence the natural projection of $w(f_1, \ldots, f_n)$ on $\prod_{i \in \Gamma} G_i$ has finite order. It follows that $f_1, \ldots, f_n$ satisfy a free word, contradicting our hypothesis. Hence $g_1, \ldots, g_n$ freely generate a free subgroup.
\end{proof}

\begin{theorem} \label{prod-main}
Let $\kappa$ be an infinite cardinal, let $\varphi : \kappa \rightarrow \Z_+$ be a function, and let $G = \prod_{i \in \kappa} S_{\varphi(i)}$ $($where $S_{\varphi(i)} = \Sym(\{1, 2, \ldots, \varphi(i)\}))$. Suppose that $\, |\{i : \varphi(i) \geq l\}| = \kappa$ for each $l \in \Z_+$. Then $G$ contains a free subgroup of rank $2^\kappa$ that is dense with respect to the product topology resulting from putting the discrete topology on each $S_{\varphi(i)}$. In particular, $G$ is almost $2^\kappa$-free.
\end{theorem}

\begin{proof}
For each $i \in \kappa$, let $\pi_i : G \rightarrow S_{\varphi(i)}$ denote the natural projection. We note that each open subset of $G$ contains a set of the form $\bigcap_{i \in I} \{f \in G : \pi_i f = g_{i}\}$, for some finite $I \subseteq \kappa$ and $g_{i} \in S_{\varphi(i)}$ ($i \in I$).

Let $\Delta$ be the set of all finite subsets of $\kappa$ (so in particular, $|\Delta| = \kappa$). For each $I \in \Delta$ let $G_I = \prod_{i \in I} S_{\varphi(i)}$, and write $G_I = \{g_{(I,j)}\}_{1\leq j \leq |G_I|}$. By Corollary~\ref{prod3}, we can find a subset $\{f_i\}_{i \in 2^\kappa} \subseteq G$ that freely generates a free group. Let us also reindex the first $k$-many $f_i$ as $\{f_i\}_{i \in \kappa} = \{f_{(I,j)} : I \in \Delta, 1\leq j \leq |G_I|\}$. Now, for each $I \in \Delta$ and $1\leq j \leq |G_I|$ define $h_{(I,j)} \in G$ so that $h_{(I,j)}$ agrees with $g_{(I,j)}$ on the coordinates indexed by $I$ and with $f_{(I,j)}$ elsewhere. Letting $H = \{h_{(I,j)} : I \in \Delta, 1\leq j \leq |G_I|\} \cup \{f_i\}_{i \in 2^\kappa \setminus \kappa},$ we see that $H$ is dense in $G$, by the observation about the open subsets of $G$ made in the previous paragraph. Also, by Lemma~\ref{prod4} (applied to all finite subsets of $H$), the elements of $H$ freely generate a free group of rank $2^\kappa$, giving us the desired conclusion.

The final claim follows from Corollary~\ref{discrete-free2}.
\end{proof}

\section{Further questions}

We conclude with a couple of questions to which we would like to know the answers.

\begin{question}
Given an integer $n > 1$, is there a topological group that is almost $n$-free but not almost $(n+1)$-free?
\end{question}

\begin{question}
Is there a completely metrizable group that is almost free but not almost $\aleph_1$-free?
\end{question}

\noindent
Department of Mathematics \newline
University of Colorado \newline
Colorado Springs, CO 80933-7150 \newline
USA \newline

\noindent Email: {\tt zmesyan@uccs.edu}


\begin{thebibliography}{00}
\bibitem{Bergman} George M.\ Bergman, \textit{Some results on embeddings of algebras, after de Bruijn and McKenzie,} Indag.\ Math.\ \textbf{18} (2007) 349--403. 

\bibitem{Bhattarcharjee} Meenaxi Bhattarcharjee, \textit{The ubiquity of free subgroups in certain inverse limits of groups,} J.\ Algebra \textbf{172} (1995) 134--146.

\bibitem{BG} E.\ Breuillard and T.\ Gelander, \textit{On dense free subgroups of Lie groups,} J.\ Algebra \textbf{261} (2003) 448--467.

\bibitem{BR} R.\ M.\ Bryant and V.\ A.\ Roman'kov, \textit{The automorphism groups of relatively free algebras,} J.\ Algebra \textbf{209} (1998) 713--723.

\bibitem{Bruijn} N.\ G.\ de Bruijn, \textit{Embedding theorems for infinite groups,} Nederl.\ Akad.\ Wetensch.\ Proc.\ Ser.\ A.\ 60 = Indag.\ Math.\ \textbf{19} (1957) 560--569.

\bibitem{Dixon} John D.\ Dixon, \textit{Most finitely generated permutation groups are free,} Bull.\ London Math.\ Soc.\ \textbf{22} (1990) 222--226.

\bibitem{Epstein} D.\ B.\ A.\ Epstein, \textit{Almost all subgroups of a Lie group are free,} J.\ Algebra \textbf{19} (1971) 261--262. 

\bibitem{Gelander} Tsachik Gelander, \textit{On deformations of $F_n$ in compact Lie groups,} Israel J.\ Math \textbf{167} (2008) 15--26.

\bibitem{GK} P.\ M.\ Gartside and R.\ W.\ Knight, \textit{Ubiquity of free subgroups,} Bull.\ London Math.\ Soc.\ \textbf{35} (2003) 624--634.

\bibitem{GMR} A.\ M.\ W.\ Glass, Stephen H.\ McCleary, and Matatyahu Rubin, \textit{Automorphism groups of countable highly homogeneous partially ordered sets,}  Math.\ Z.\ \textbf{214} (1993) 55--66. 

\bibitem{MS} Garvin Melles and Saharon Shelah, \textit{$\Aut(M)$ has a large dense free subgroup for saturated $M$,} Bull.\ London Math.\ Soc.\ \textbf{26} (1994) 339--344.

\end{thebibliography}
\end{document}